\documentclass{amsart}
 
\usepackage{amsmath,amssymb}

\newtheorem{theorem}{Theorem}
\newcommand{\bt}{\begin{theorem}}
\newcommand{\et}{\end{theorem}}

\newtheorem{lemma}{Lemma}
\newcommand{\bl}{\begin{lemma}}
\newcommand{\el}{\end{lemma}}

\newtheorem{corollary}{Corollary}
\newcommand{\bc}{\begin{corollary}}
\newcommand{\ec}{\end{corollary}}

\newtheorem{example}{Example}
\newcommand{\bex}{\begin{example}}
\newcommand{\eex}{\end{example}}

\newtheorem{problem}{Problem} 
\newcommand{\bprob}{\begin{problem}}
\newcommand{\eprob}{\end{problem}}

\newtheorem{definition}{Definition}
\newcommand{\bdf}{\begin{definition}}
\newcommand{\edf}{\end{definition}}

\newcommand{\bpf}{\begin{proof}}
\newcommand{\epf}{\end{proof}}
\newcommand{\beq}{\begin{equation}}
\newcommand{\eeq}{\end{equation}}
\newcommand{\benum}{\begin{enumerate}}
\newcommand{\eenum}{\end{enumerate}}

\newcommand{\N}{\ensuremath{ \mathbf N }}
\newcommand{\NO}{\ensuremath{ {\mathbf N}_0 }}

\newcommand{\X}{\ensuremath{\mathbf X}}
\newcommand{\Z}{\ensuremath{\mathbf Z}}

\newcommand{\card}{\text{card}}
\newcommand{\bq}{\begin{eqnarray*}}
\newcommand{\eq}{\end{eqnarray*}}

\title[Additive Representation Functions]{Inverse Problems for Representation Functions in Additive Number Theory}

\date{\today}



\author{Melvyn B. Nathanson}
\thanks{The work of M.B.N. was supported
in part by grants from the NSA Mathematical Sciences Program and
the PSC-CUNY Research Award Program}
\address{Department of Mathematics,
Lehman College (CUNY),
Bronx, New York 10468}
\email{melvyn.nathanson@lehman.cuny.edu}
\curraddr{School of Mathematics, Institute for Advanced Study,
Princeton, NJ 08540}
\email{melvyn@ias.edu}

\begin{document}

\begin{abstract}
For every positive integer $h$, the representation function of order $h$ associated to a subset $A$ of the integers or, more generally, of any group or semigroup $\X$, counts the number of ways an element of $\X$ can be written as the sum (or product, if $\X$ is nonabelian) of $h$ not necessarily distinct elements of $\X$.  The direct problem for representation functions in additive number theory begins with a subset $A$ of $\X$ and seeks to understand its representation functions.  The inverse problem for representation functions starts with a function $f:\X \rightarrow \N_0 \cup \{\infty\}$ and asks if there is a set $A$ whose representation function is $f$, and, if the answer is yes, to classify all such sets.  This paper is a survey of recent progress on the inverse representation problem.  
\end{abstract}

\maketitle

\section{Asymptotic density}

Let $\N, \N_0,$ and \Z\ denote, respectively, the sets of positive integers, nonnegative integers, and all integers.

For any set $A$ of integers, we define the \emph{counting function} $A(y,x)$ of $A$ by
\[
A(y,x) = \sum_{\substack{a\in A\\ y \leq a \leq x}} 1
\]
for all real numbers $x$ and $y$.
We define the \emph{upper asymptotic density}
\[
d_U(A) = \limsup_{x\rightarrow\infty} \frac{A(-x,x)}{2x+1}
\]
and the \emph{lower asymptotic density}
\[
d_L(A) = \liminf_{x\rightarrow\infty} \frac{A(-x,x)}{2x+1}.
\]
The set $A$ has \emph{asymptotic density} $d(A) = \alpha$ if 
$d_U(A) = d_L(A) = \alpha$ or, equivalently, if
\[
d(A) = \lim_{x\rightarrow\infty} \frac{A(-x,x)}{2x+1} = \alpha.
\]
Let $B = \Z \setminus A.$  Then $d_U(A) = \alpha$ if and only if $d_L(A) = 1 - \alpha.$
If $S$ and $W$ are sets of integers, then the set $W$ has \emph{relative upper asymptotic density $d_U(W,S) = \alpha$ with respect to $S$} if 
\[
\limsup_{x\rightarrow\infty} \frac{(W \cap S)(-x,x)}{S(-x,x)} = \alpha.
\]
Relative lower asymptotic density $d_L(W,S)$ and relative density $d(W,S)$ are defined similarly.  In particular, if $S= \N$ and $A$ is a set of positive integers, then $A$ has relative asymptotic density $\alpha$ with respect to \N\ if 
\[
\lim_{x\rightarrow\infty} \frac{A(-x,x)}{\N(-x,x)} = 
\lim_{x\rightarrow\infty} \frac{A(1,x)}{[x]} = \alpha.
\]

\section{Sumsets and bases}
Let $A_1$ and $A_2$ be subsets of an additive abelian semigroup $\X$.  We define the \emph{sumset} 
\[
A_1 + A_2 = \{a_1 + a_2 : a_1\in A_1 \text{ and } a_2 \in A_2\}.
\]
For every positive integer $h \geq 3$, if $A_1, A_2,\ldots, A_h$ are subsets of \X, then we define the sumset $A_1 + \cdots +A_{h-1} +  A_h$ inductively by
\[
A_1 + \cdots +A_{h-1} +  A_h = (A_1 + \cdots +A_{h-1}) +  A_h.
\]
If $A = A_i$ for $i=1,\ldots, h,$ then we write
\[
hA = \underbrace{A+\cdots + A}_{\text{$h$ times}}.
\]
The set $hA$ is called the \emph{$h$-fold sumset} of $A$.

We define $0A = \{ 0 \}.$

If $A \subseteq \X$ and $x\in\X,$ we define the \emph{shift} 
$A+x = A+\{ x \}.$

A central concept in additive number theory is \emph{basis}.  Let $S$ be a subset of \X.  The set $A$ is called 
\benum
\item
a \emph{basis of order $h$ for $S$} if $S \subseteq hA,$ that is, if every element of $S$ can be represented as the sum of $h$ not necessarily distinct elements of $A$,
\item
an \emph{asymptotic basis of order $h$ for $S$} if $S\setminus hA$ is finite, that is, if all but finitely many elements of $S$ can be represented as the sum of $h$ not necessarily distinct elements of $A$.
\eenum
For example,  the nonnegative cubes are a basis of order 9 for $\N_0$ (Wieferich's theorem), an asymptotic basis of order 7 for $\N_0$ (Linnik's theorem), and a basis of order 4 for almost all $\N_0$ (Davenport's theorem).   A large part of classical additive number theory is the study of how special sets of integers (for example, the $k$-th powers, polygonal numbers, and primes) are bases for the nonnegative integers (cf. Nathanson~\cite{nath96aa}).  

Our definition of basis is weak in the sense that, if  \X\ is an abelian semigroup with additive identity 0, then every subset of \X\ has a basis of order $h$ for all $h \geq 1.$.   The reason is that $0 \in \X$ implies that $h\X = \X$ for every positive integer $h$, and so, if $S$ is any subset of \X, then $S \subseteq h\X$. 
We shall call the subset $A$ of \X\ 
\benum
\item
an \emph{exact basis of order $h$ for $S$} if $S = hA,$ that is, if the  elements of $S$ are precisely the elements of \X\ that can be represented as the sum of $h$ not necessarily distinct elements of $A$,
\item
an \emph{exact asymptotic basis of order $h$ for $S$} if $hA \subseteq S$ and $S\setminus hA$ is finite.
\eenum
In additive subsemigroups of the integers, the set $A$ is a \emph{basis of order $h$ for almost all $S$} if $S\setminus hA$ has relative asymptotic density zero with respect to $S$, and an \emph{exact basis of order $h$ for almost all $S$} if $hA \subseteq S$ and $S\setminus hA$ has relative asymptotic density zero with respect to $S$.

\section{Direct and inverse problems for sumsets}
Let \X\ be an additive abelian semigroup.  Given subsets $A_1, \ldots, A_h$ of \X, a \emph{direct problem} in additive number theory is to describe the sumset $A_1 + \cdots + A_h.$  In particular, for any $A \subseteq \X,$ the direct problem is to describe the $h$-fold sumsets $hA$ for all $h \geq 2.$  If \X\ contains an additive identity 0 and if $0 \in A\subseteq \X,$ then we obtain an increasing sequence of sumsets 
\beq   \label{RS:hAsequence}
A \subseteq 2A \subseteq  \cdots \subseteq hA \subseteq (h+1)A \subseteq \cdots.
\eeq
An important open problem is to describe the evolution of structure in the sequence $\{hA\}_{h=1}^{\infty}$.  For example, let $\X = \N_0$ be the additive semigroup of nonnegative integers.  Let $A$ be a set of nonnegative integers such that $d_L(h_0A) > 0$ for some positive integer $h_0.$   By translation and contraction, we can assume that $0 \in A$ and $\gcd(A)=1.$  The the sequence~\eqref{RS:hAsequence} eventually stabilizes as a co-finite subset of \NO, that is, there exists an integer $h_1 \geq h_0$ such that $h_1A$ contains all sufficiently large integers and $hA = h_1A$ for all $h \geq h_1$  (Nash-Nathanson~\cite{nath85a}).   However, if $d_L(hA)=0$ for all positive integers $h$, then the structure of the sumsets $hA$ is mysterious.  It must happen that very regular infinite configurations of integers develop in the sumsets, but nothing is known about them.  

The simplest \emph{inverse problem for sumsets} is:\\
\begin{center}
\framebox{ \Large{What sets are sumsets? } }
\end{center}
\vspace{0.3cm}
This can be called the \emph{sumset recognition problem}:  Given a subset $S$ of the abelian semigroup \X\ and an integer $h \geq 2$, do there exist subsets $A_1,\ldots, A_h \subseteq \X$ such 
\[
S = A_1 + \cdots + A_h?
\]
Similarly, we have \emph{basis recognition problems}.  Let $S$ be a subset of the abelian semigroup \X, and let $h\geq 2.$
Does there exist an exact basis of order $h$ for $S$, that is, a set $A \subseteq \X$ such that $hA = S$?  If $S$ does have an exact basis, describe the set
\[
\mathcal{E}_h(S) = \{A \subseteq \X: S = hA\}.
\]
More generally, do there exist exact asymptotic bases of order $h$ for $S$?  If so, describe the set
\[
\mathcal{E}_h^{\text{asy}}(S) = \{A \subseteq \X: hA \subseteq S \text{ and } \card(S\setminus hA) < \infty\}.
\]

\section{Representation functions of semigroups}
Let $\X$ be an abelian semigroup, written additively.  For $A\subseteq \X$, let $A^h$ denote the set of all $h$-tuples of $A.$   Two $h$-tuples $(a_1,\ldots,a_h)\in \X^h$ and $(a'_1,\ldots,a'_h)\in \X^h$ are equivalent if there is a permutation $\tau:\{1,\ldots,h\} \rightarrow \{1,\ldots,h\}$ such that $\tau(a_i) = a'_i$ for $i=1,\ldots,h.$  If $\X$ is the semigroup of integers or nonnegative integers (or if $\X$ is any totally ordered set), then every equivalence class contains a unique $h$-tuple $(a_1,\ldots,a_h)$ such that $a_i \leq a_{i+1}$ for $i=1,\ldots,h-1.$

Let $A_1,\ldots,A_h$ be subsets of $\X$ and let $x$ be an element of $\X.$   We define the \emph{ ordered representation function}
\[
R_{A_1,\ldots,A_h}(x) = \card\left( \left\{ (a_1,\ldots,a_h) \in A_1 \times \cdots \times A_h : a_1 + \cdots + a_h = x \right\} \right).
\]
If $A_i = A$ for $i=1,\ldots,h,$ then we write 
\[
R_{A,h}(x)  = \card\left( \left\{ (a_1,\ldots,a_h) \in A^h : a_1 + \cdots + a_h = x \right\} \right).
\]
Two other representation functions arise often and naturally in additive number theory. 
The \emph{unordered representation function} $r_{A,h}(x)$ counts the number of equivalence classes of $h$-tuples $(a_1,\ldots,a_h)$ such that $a_1+\cdots + a_h = x.$ 
The \emph{unordered restricted representation function}\footnote{We could also introduce an \emph{ordered restricted representation function} $\hat{R}_{A,h}(x)$ that counts the number of $h$-tuples $(a_1,\ldots,a_h)$ of pairwise distinct elements of $\X$ such that $a_1+\cdots + a_h = x.$  This is unnecessary, however, because $\hat{R}_{A,h}(x) = h!\hat{r}_{A,h}(x)$ for all $x \in \X.$  The relation between the ordered and unordered representation functions $R_{A,h}(x)$ and $r_{A,h}(x)$ is more complex.}
$\hat{r}_{A,h}(x)$ counts the number of equivalence classes of $h$-tuples $(a_1,\ldots,a_h)$ of pairwise distinct elements of $\X$ such that $a_1+\cdots + a_h = x.$ 

If \X\ is a subsemigroup of the integers or of any totally ordered semigroup, then 
\[
r_{A,h}(x)  = \card\left( \left\{ (a_1,\ldots,a_h) \in A^h : a_1 \leq \cdots \leq a_h \text{ and } a_1 + \cdots + a_h = x \right\} \right)
\]
and
\[
\hat{r}_{A,h}(x)  = \card\left( \left\{ (a_1,\ldots,a_h) \in A^h : a_1 <\cdots < a_h \text{ and } a_1 + \cdots + a_h = x \right\} \right).
\]

\section{Direct and inverse problems for representation functions}  
A fundamental  \emph{direct problem} in additive number theory is to describe the representation functions of finite and infinite subsets of the integers and of other abelian semigroups.    For example, if $\X = \N_0 = A,$ then  $R_{A,2}(n) = n+1$ and $r_{A,2}(n) = [(n+2)/2]$  for all $n\in \N_0.$  If $\X = \Z = A,$ then $R_{A,2}(n)  =r_{A,2}(n) = \infty$ for all $n\in \Z.$  
More generally, we ask:  Given a semigroup $\X$, a family $\mathcal{A}$ of subsets of $\X$, and a positive integer $h,$ what properties are shared by all of the representation functions  associated with sets $A \in \mathcal{A}$?  
These are direct problems.

The simplest \emph{inverse problem for representation functions} is:
\vspace{0.3cm}
\begin{center}
\framebox{ \Large{What functions are representation functions? } }
\end{center}
\vspace{0.3cm}
More precisely, if $\mathcal{A}$ is a family of subsets of a semigroup $\X$ and if $h$ is a positive integer, let $\mathcal{R}_h^{\text{ord}}(\mathcal{A})$ and $\mathcal{R}_h^{\text{unord}}(\mathcal{A})$ denote, respectively, the sets of ordered and ordered representation functions of order $h$ associated with sets $A \in \mathcal{A},$ that is, 
\[
\mathcal{R}^{\text{ord}}_h(\mathcal{A}) = \{R_{A,h} : A \in \mathcal{A} \}
\]
and
\[
\mathcal{R}^{\text{unord}}_h(\mathcal{A}) = \{r_{A,h} : A \in \mathcal{A} \}.
\]
The inverse problem is to determine if a given function $f$ is a representation function, and, if so, to describe all sets $A \in \mathcal{A}$ such that $R_{A,h} = f$  or $r_{A,h} = f.$  
This is particularly interesting when $\mathcal{A}$ is the set of bases or asymptotic bases of order $h$ for $\X$.

There is an important difference between the representation functions of asymptotic bases for the integers and the nonnegative integers.  If $f$ is the unordered representation function of an asymptotic basis for a semigroup \X, then the set $f^{-1}(0)$ is finite.   If $\X = \Z,$ then a fundamental theorem in additive number theory (Theorem~\ref{RS:theorem:FundRep}) states that for every $h \geq 2$ and for every function $f:\Z \rightarrow \NO \cup \{\infty\}$ with $\card(f^{-1}(0)) < \infty$, there exists a set $A$ such that $r_{A,h}(n) = f(n)$ for all $n\in \Z.$  Equivalently, if $\mathcal{A}$ is the set of all asymptotic bases of order $h$ for \Z, then 
\[
\mathcal{R}^{\text{unord}}_h(\mathcal{A}) = \{ f:\Z \rightarrow \N_0 \cup \{\infty\} : \card(f^{-1}(0)) < \infty \}.
\]
For the semigroup of nonnegative integers, however, it is false that every function $f:\NO \rightarrow \NO$ with only finitely many zeros is the unordered representation function for an asymptotic basis of order $h$.   Indeed, very little is know about representation functions of asymptotic bases of finite order for \NO.

\section{Representation functions for sets of nonnegative integers}
If $A$ is a set of nonnegative integers, then for every positive integer $h$  the number of representations of an integer as the sum of $h$ elements of $A$ is finite.   We introduce the following three sets of arithmetic functions:
\[
\mathcal{F}(\NO) = \left\{f:\NO \rightarrow \NO  \right\}
\]
\[
\mathcal{F}_{\infty}(\NO) = \left\{f:\NO \rightarrow \NO : \text{$f^{-1}(0)$ is a set of density 0}\right\}
\]
and
\[
\mathcal{F}_0(\NO) = \left\{f:\NO \rightarrow \NO : \text{$f^{-1}(0)$ is a finite set}\right\}.
\]
Then
\[
\mathcal{F}_0(\NO) \subset \mathcal{F}_{\infty}(\NO)  \subset \mathcal{F}(\NO).
\]
For $h \geq 2,$  the set $\mathcal{F}_0(\NO)$ contains the representation functions of all bases and asymptotic bases of order $h$ for \NO, and the  set $\mathcal{F}_{\infty}(\NO)$ contains the representation functions of all bases of order $h$ for almost all \NO.

\bprob
Let $h \geq 2.$  Find necessary and sufficient conditions for a function in $\mathcal{F}_0$ to be the representation function for an asymptotic basis of order $h$ for \NO.
\eprob

\bprob
Let $h \geq 2.$  Find necessary and sufficient conditions for a function in $\mathcal{F}_{\infty}$ to be the representation function for a basis of order $h$ for almost all \NO.
\eprob

\bprob
Let $h \geq 2.$  Find necessary and sufficient conditions for a function in $\mathcal{F}$ to be the representation function for a subset of \NO.
\eprob

We can also count the number of representations of a nonnegative integer as the sum of a bounded number of elements of a set that contains both nonnegative and negative integers.

\bprob
Let $h \geq 2.$  Find necessary and sufficient conditions for a function in $\mathcal{F}$ to be the representation function for the nonnegative integers in the $h$-fold sumset of  a subset of \Z.
\eprob

We can express the ordered and unordered representation functions of a set of nonnegative integers in terms of generating functions.
Define the \emph{generating function} for the set $A$ of nonnegative integers as the power series 
\[
G_A(z)=\sum_{a\in A} z^a.
\]
This  can be used both as a formal power series and as an analytic function that converges for $|z|<1.$ 
We have the identities
\[
\sum_{n=0}^{\infty} r_{A,2}(n)z^n = \frac{1}{2}\left( G^2_A(z) + G_A(z^2) \right)
\]
\[
\sum_{n=0}^{\infty} \hat{r}_{A,2}(n)z^n = \frac{1}{2}\left( G^2_A(z) - G_A(z^2) \right)
\]
and, for all $h \geq 1,$
\[
\sum_{n=0}^{\infty} R_{A,h}(n)z^n =  G^h_A(z).
\]

If $A$ is a set of integers, then the ordered representation function $R_{A,2}(n)$ is odd if and only if $n$ is even and $n/2\in A.$  It follows that  $R_{A,2}(n)$ is eventually constant only if and only if $A$ is finite.  Moreover, the ordered representation function $R_{A,2}$ uniquely determines the set $A$.  Thus, for every function $f \in \mathcal{F}(\NO),$ there exists at most one set $A$ such that $R_{A,2} = f.$  Theorem~\ref{RS:theorem:nath-unique}  generalizes this observation to all $h \geq 2.$

It is also true that the unordered representation function $r_{A,2}(n)$ for a set $A$ of nonnegative integers is eventually constant only if $A$ is finite. 

\bt[Dirac~\cite{dira51}]
If $\mathcal{A}$ is an infinite set of  nonnegative integers, then the representation function  $r_{A,2}(n)$ is not eventually constant.
\et

\bpf
Let $A$ be an infinite set of nonnegative integers such that $r_{A,2}(n) = c $ for all $n \geq n_0.$  Since $A$ is infinite, we have $r_{A,2}(2a) \geq 1$ for all $a\in A,$ and so $c \geq 1.$  There is a polynomial $P(z)$ such that
\begin{align*}
\frac{1}{2}\left( G^2_A(z) + G_A(z^2) \right)  
& = \sum_{n=0}^{\infty} r_{A,2}(n)z^n \\
& = \sum_{n=0}^{n_0-1} r_{A,2}(n)z^n + \sum_{n=n_0}^{\infty} cz^n \\
& = \frac{P(z)}{1-z}.
\end{align*}
  Let $0 < x < 1$ and $z=-x.$  Then $G_A(z) = G_A(-x)$ is real and so $G^2_A(z)  \geq 0$ and 
\[
 \frac{2P(-x)}{1+x} =  G^2_A(-x) + G_A(x^2) \geq G_A(x^2).
 \] 
Taking the limit as $x\rightarrow 1^-,$ we see that the left side of this equality converges to $P(-1)$ but the right side diverges to infinity.  This is impossible, and so the representation function $r_{A,2}(n)$ cannot be eventually constant.  
\epf

Dirac's theorem is a special case of a famous unsolved problem in additive number theory.   Erd\H os and Tur\' an~\cite{erdo-tura41} conjectured that if $A$ is an asymptotic basis of order 2 for the  nonnegative integers, then $\limsup_{n\rightarrow\infty} r_{A,2}(n) =\infty.$  This conjecture is itself only a small part of the problem of characterizing the representation functions of additive bases of finite order for \NO.  It is interesting to note that the modular analogue of the Erd\" os-Tur\' an  conjecture is false.

\bt[Tang-Chen~\cite{tang-chen06}]
There is an integer $m_0$ such that, for every $m \geq m_0,$ there is a set $A_m \subseteq \Z/m\Z$ such that $A_m$ is a basis of order 2 for $\Z/m\Z$ and $r_{A_m,2}(x) \leq 768$ for all $x \in \Z/m\Z.$
\et

It is also interesting that the multiplicative Erd\H os and Tur\' an conjecture is true.  If  $A$ is a set of positive integers such that every sufficiently large positive integer is the product of two elements of $A$, then the number of representations of an integer $n$ as the product of two elements of $A$ is unbounded (Erd\H os~\cite{erdo64a}, Ne{\u s}et{\u r}il and R{\" o}dl~\cite{nese-rodl85}, Nathanson~\cite{nath87c}).

\subsection{Ordered representation functions}
The first inverse theorems for ordered representation functions of sets of nonnegative integers are the following.   

\bt[Nathanson~\cite{nath78d}]    \label{RS:theorem:nath-unique}
Let $h \geq 2.$  
If $A$ and $B$ are sets of nonnegative integers such that $R_{A,h}(n)=R_{B,h}(n)$ for all $n\in \N_0,$  then $A=B.$
\et

\begin{proof}
Since $A=\emptyset$ if and only if $B=\emptyset,$ we can assume that both $A$ and $B$ are nonempty sets.   Then the generating functions
\[
G_A(z)=\sum_{a\in A} z^a  \qquad\text{and} \qquad G_B(z)=\sum_{b\in B} z^b
\]
are nonzero power series with nonnegative coefficients.   We have 
\[
G^h_A(z)  = \left(\sum_{a\in A} z^a  \right)^h  = \sum_{n=0}^{\infty} R_{A,h}(n) = \sum_{n=0}^{\infty} R_{B,h}(n)  = \left(\sum_{b\in B} z^b  \right)^h  = G^h_B(z)
\]
and so
\[
0 = G^h_A(z) - G^h_B(z) = \left( G_A(z) - G_B(z) \right)\left( \sum_{i=0}^{h-1} G^{h-1-i}_A(z)  G_B^i(z) \right).
\]
The coefficients of the power series $\sum_{i=0}^{h-1} G^{h-1-i}_A(z)  G_B^i(z)$ are nonnegative and not all zero, hence this series is nonzero and so $ G_A(z) - G_B(z) = 0.$  This implies that $A=B.$
\end{proof}

Let $A^{\ast},B^{\ast}$, and $T$ be finite sets of integers.  If each residue class modulo $m$ contains exactly the same number of elements of $A^{\ast}$ as elements of $B^{\ast}$, then we write $A^{\ast} \equiv B^{\ast} \pmod{m}.$  If for each integer $n$  the number of pairs $(a,t)  \in A^{\ast}\times  T$ such that  $a+t \equiv n \pmod{m}$ equals the number of 
 pairs $(b,t)  \in B^{\ast}\times  T$ such that $b+t \equiv n \pmod{m}$, then we write
\[
A^{\ast} +T \equiv B^{\ast} +T \pmod{m}.
\]

\bt[Nathanson~\cite{nath78d}]  \label{RS:theorem:NathOrder}
Let $A$ and $B$ be sets of nonnegative integers.  Then $R_{A,2}(n) = R_{B,2}(n)$ for all sufficiently large $n$ if and only if there exist
\benum
\item[(i)]
a nonnegative integer $n_0$ and sets $A^{\ast},B^{\ast} \subseteq \{ 0,1,2,\ldots,n_0\}$, and
\item[(ii)]
a positive integer $m$ and a set $T \subseteq \{ 0,1,2,\ldots,m-1\}$ with
\[
A^{\ast} +T \equiv B^{\ast} +T \pmod{m}
\]
\eenum
such that 
\beq  \label{RS:AB}
A = A^{\ast} \cup C \qquad\text{and}\qquad B = B^{\ast} \cup C
\eeq
where
\beq  \label{RS:C}
C = \{ c \in \NO : c > n_0 \text{ and } c \equiv t\pmod{m} \text{ for some $t\in T$}\}.
\eeq
\et

\begin{proof}
Let $n_0$ and $m$ be integers and let $A^{\ast}, B^{\ast},$ and $T$ be finite sets of integers satisfying conditions~(i) and~(ii).  Define the sets $A,B,$ and $C$ by~\eqref{RS:AB} and~\eqref{RS:C}.  Since $A^{\ast} \cap C = \emptyset$ and $B^{\ast} \cap C = \emptyset$, it follows that for every integer $n$ we have
\[
R_{A,2}(n) = R_{A^{\ast},2}(n) + 2R_{A^{\ast},C}(n) + R_{C,2}(n)
\]
and
\[
R_{B,2}(n) = R_{B^{\ast},2}(n) + 2R_{B^{\ast},C}(n) + R_{C,2}(n)
\]
where $R_{A^{\ast},C}(n)$ (resp. $R_{B^{\ast},C}(n)$ ) is the number of ordered pairs $(a^{\ast},c)\in A^{\ast} \times C$ (resp. $(b^{\ast},c)\in B^{\ast} \times C$) such that $a^{\ast}+c =n$ (resp. $b^{\ast}+c=n$).

Let $n>2n_0.$  Since $\max\left( A^{\ast} \cup B^{\ast} \right) \leq n_0,$ it follows that  $R_{A^{\ast},2}(n) = R_{B^{\ast},2}(n) = 0$ and so
 $R_{A,2}(n)= R_{B,2}(n)$ if and only if $R_{A^{\ast},C}(n) = R_{B^{\ast},C}(n).$
If $a^{\ast}\in A,$ then $n-a^{\ast} > 2n_0 - a^{\ast} \geq n_0.$  
It follows that $n-a^{\ast}\in C$ if and only if $n-a^{\ast} \equiv t \pmod{m}$ for some $t\in T$.  Since $A^{\ast} +T \equiv B^{\ast} +T \pmod{m},$ it follows that
\begin{align*}
R_{A^{\ast},C}(n) 
& = \card\left( \left\{ (a^{\ast},c) \in A^{\ast} \times C : a^{\ast} + c = n \right\}\right) \\
& = \sum_{t\in T}\card\left( 
\left\{ ( a^{\ast},c) \in A^{\ast} \times C : a^{\ast} + c = n 
\text{ and } c\equiv t \pmod{m}\right\} \right) \\
& = \sum_{t\in T}\card\left( 
\left\{ ( a^{\ast},t) \in A^{\ast} \times T : a^{\ast} + t \equiv n \pmod{m}\right\} \right) \\
& = \sum_{t\in T}\card\left( 
\left\{ ( b^{\ast},t) \in B^{\ast} \times T : b^{\ast} + t \equiv n \pmod{m}\right\} \right) \\
& = \sum_{t\in T}\card\left( 
\left\{ ( b^{\ast},c) \in B^{\ast} \times C : b^{\ast} + c = n 
\text{ and } c\equiv t \pmod{m}\right\} \right) \\
& = R_{B^{\ast},C}(n)
\end{align*}
Thus,  the representation functions of the sets $A$ and $B$ eventually coincide.

Conversely, let $A$ and $B$ be distinct sets of integers such that $R_{A,2}(n) = R_{B,2}(n)$ for all integers $n > n_1$.  Since $A$ is finite if and only if $R_{A,2}(n) = 0$ for all sufficiently large $n,$ it follows that the representation functions of any pair of finite sets eventually coincide, and so $A$ is finite if and only if $B$ is finite.  Thus, we can set $A^{\ast} = A$, $B^{\ast} = B$, and  $T = C=\emptyset$.

Suppose that $A$ and $B$ are distinct infinite sets of integers.   Applying the generating functions $G_A(z) = \sum_{a\in A} z^a$ and $G_B(z) = \sum_{b \in B} z^b$, we have
\[
G^2_A(z) - G^2_B(z) = \sum_{n=0}^{\infty} \left( R_{A,2}(n) - R_{B,2}(n) \right) z^n = P(z)
\]
where $P(z)$ is a polynomial of degree at most $n_1$.  
The ordered representation function $R_{A,2}(n)$ (resp. $R_{B,2}(n)$) is odd if and only if $n$ is even and $n/2 \in A$ (resp. $n/2 \in B$).    It follows that the sets $A$ and $B$ coincide for $n > n_1/2$, and so there is a nonzero polynomial $Q(z)$ of degree at most $n_1/2$ such that 
\[
G_A(z) - G_B(z) = Q(z).
\] 
We obtain a rational function
\[
 G_A(z) +  G_B(z) = \frac{G^2_A(z) - G^2_B(z)}{G_A(z) - G_B(z)} = \frac{P(z)}{Q(z)}.
\]
Therefore, the coefficients of the power series $G_A(z)+G_B(z)$ satisfy a linear recurrence relation.  For $n > n_1/2,$  the coefficient of $z^n$ in $G_A(z) + G_B(z)$ is 2 if $n \in A \cap B$ and 0 if $n \notin A\cap B$.   Since a sequence defined by a linear recurrence in a finite set must be eventually periodic, it follows that there are positive integers $m$ and $n_0$ and a set $T \subseteq \{0,1,\ldots, m-1 \}$ such that, for $n > n_0,$ we have $n \in A\cap B$ if and only if $n \equiv t \pmod{m}$ for some $t\in T.$   Let
\[
C = \{ c \in \NO : c > n_0 \text{ and } c \equiv t\pmod{m} \text{ for some $t \in T$} \}.
\]
Let $A^{\ast} = A \cap [0,n_0]$ and $B^{\ast} = B \cap [0,n_0]$.  Then $A^{\ast} \cap C  = B^{\ast} \cap C = \emptyset,$ and 
$A = A^{\ast} \cup C$ and $B = B^{\ast} \cup C.$
For $n > 2n_0$ we have
\[
2R_{A^{\ast},C}(n)  = R_{A,2}(n) - R_{C,2}(n) =  R_{B,2}(n) - R_{C,2}(n) =  2R_{B^{\ast},C}(n)
\]
where, as above, $R_{A^{\ast},C}(n)$ (resp. $R_{B^{\ast},C}(n)$) is the number of solutions of the congruence $n \equiv a + t \pmod{m}$ (resp. $n \equiv b + t \pmod{m}$) with $t\in T$ and $a\in A^{\ast}$ (resp. $b \in B^{\ast}$).  Therefore, $A^{\ast} +T \equiv B^{\ast} +T \pmod{m},$ and the Theorem follows.
\end{proof}

\bprob
Let $h \geq 3.$  Describe all pairs of sets of nonnegative integers whose ordered representation functions of order $h$ eventually coincide.  Equivalently, classify all pairs 
$(A,B)$ of sets of nonnegative integers such that $R_{A,h}(n) = R_{B,h}(n)$ for all sufficiently large integers $n$.  
\eprob

\subsection{Unordered representation functions}
Theorem~\ref{RS:theorem:NathOrder} completely describes all pairs of sets of nonnegative integers whose ordered representation functions of order 2 eventually coincide.  The analogous problem for unordered representation functions is open.

\bprob
Describe all pairs of sets of nonnegative integers whose unordered representation functions of order 2 eventually coincide. \eprob

\bprob
Let $h \geq 3.$  Describe all pairs of sets of nonnegative integers whose unordered representation functions eventually coincide.  \eprob

The behavior of unordered representation functions is more exotic than that of ordered representation functions.  For example, the following beautiful result describes partitions of the nonnegative integers into disjoint sets $A$ and $B$ whose unordered representation functions eventually coincide.  

\bt[S\' andor~\cite{sand04}]
Let $A$ be a set of nonnegative integers, and let $B = \N_0 \setminus A.$  There exists a positive integer $N$ such that  $r_{A,2}(n) = r_{B,2}(n)$ for all $n \geq 2N-1$ if and only if 
\benum
\item[(i)]
\[
\card( A \cap [0,2N-1] ) = N
\]
\item[(ii)]
for every integer $a \geq N$,
\[
a \in A \text{ if and only if } 2a \notin A
\]
and
\[
a \in A \text{ if and only if } 2a +1\in A
\]
\eenum
\et

\begin{proof}
Let $\chi_A(n)$ denote the characteristic function of the set $A$, that is,
\[
\chi_A(n) = 
\begin{cases}
1 & \text{if $n\in A$} \\
0 & \text{if $n\notin A$.} 
\end{cases}
\]
Since $B = \N_0 \setminus A,$ we have 
\[
\chi_B(n) = 1 - \chi_A(n) \qquad\text{for all $n\in \N_0$.}
\]
Defining the generating functions
\[
G_A(z) = \sum_{a\in A} z^a = \sum_{n=0}^{\infty} \chi_A(n) z^n
\]
and
\[
G_B(z) = \sum_{b\in B} z^b 
= \sum_{n=0}^{\infty}\left( 1 - \chi_A(n) \right)z^n
= \frac{1}{1-z} - G_A(z)
\]
we obtain
\[
\sum_{n=0}^{\infty} r_{A,2}(n)z^n = \frac{1}{2}\left(G_A(z)^2 + G_A(z^2)\right)
\]
and
\begin{align*}
\sum_{n=0}^{\infty} r_{B,2}(n)z^n 
& = \frac{1}{2}\left(G_B(z)^2 + G_B(z^2)\right) \\
& = \frac{1}{2}\left( \left( \frac{1}{1-z} - G_A(z)\right)^2 + \left(\frac{1}{1-z^2} - G_A(z^2) \right)\right)  \\
& = \frac{1}{2}\left( \frac{2}{(1-z^2)(1-z)} - \frac{2G_A(z)}{1-z} + G_A(z)^2 - G_A(z^2) \right) \\
& = \frac{1}{2}\left(G_A(z)^2 + G_A(z^2)\right)
+\left(  \frac{1}{(1-z^2)(1-z)} - \frac{G_A(z)}{1-z} - G_A(z^2)  \right) \\
& = \sum_{n=0}^{\infty} r_{A,2}(n)z^n 
+ \frac{1}{1-z}\left(  \frac{1}{1-z^2} - G_A(z) - (1-z)G_A(z^2)  \right) \\
& = \sum_{n=0}^{\infty} r_{A,2}(n)z^n 
+ \frac{1}{1-z}\left( \sum_{n=0}^{\infty}  z^{2n} 
 - \sum_{n=0}^{\infty}  \chi_A(n)z^{n}   - \sum_{n=0}^{\infty}  \chi_A(n)z^{2n}  +  \sum_{n=0}^{\infty}  \chi_A(n)z^{2n+1}  \right)\\
& = \sum_{n=0}^{\infty} r_{A,2}(n)z^n 
+ \frac{1}{1-z}\left( \sum_{n=0}^{\infty} \left(1 - \chi_A(n)- \chi_A(2n)\right) z^{2n} + \sum_{n=0}^{\infty} \left(\chi_A(n)- \chi_A(2n+1)\right) z^{2n+1} \right).
\end{align*}
We define the function
\[
Q(z) = \sum_{n=0}^{\infty}\left(  r_{A,2}(n) - r_{B,2}(n)\right) z^n.
\]
Then
\begin{align*}
(1-z)Q(z) =  
&  \sum_{n=0}^{\infty} \left(1 - \chi_A(n)- \chi_A(2n)\right) z^{2n} + \sum_{n=0}^{\infty} \left(\chi_A(n)- \chi_A(2n+1)\right) z^{2n+1} \\
  =  & \sum_{n=0}^{N-1} \left(1 - \chi_A(n)- \chi_A(2n)\right) z^{2n} + \sum_{n=0}^{N-1} \left(\chi_A(n)- \chi_A(2n+1)\right) z^{2n+1} + \\
& +  \sum_{n=N}^{\infty} \left(1 - \chi_A(n)- \chi_A(2n)\right) z^{2n} + \sum_{n=N}^{\infty} \left(\chi_A(n)- \chi_A(2n+1)\right) z^{2n+1}.
 \end{align*}
Let $N$ be a positive integer.  We have $r_{A,2}(n) = r_{B,2}(n)$ for all $n \geq 2N-1$ if and only if $Q(z)$ is a polynomial  of degree at most $2N-2$.   Then $(1-z)Q(z)$ has degree at most $2N-1,$ and we have the two equations
\[
(1-z)Q(z) =  \sum_{n=0}^{N-1} \left(1 - \chi_A(n)- \chi_A(2n)\right) z^{2n} + \sum_{n=0}^{N-1} \left(\chi_A(n)- \chi_A(2n+1)\right) z^{2n+1} 
\]
and
\[
0 =  \sum_{n=N}^{\infty} \left(1 - \chi_A(n)- \chi_A(2n)\right) z^{2n} + \sum_{n=N}^{\infty} \left(\chi_A(n)- \chi_A(2n+1)\right) z^{2n+1}.
\]
If the first equation holds, then, setting $z=1$, we obtain 
\[
0  =  \sum_{n=0}^{N-1} \left(1 -  \chi_A(2n) - \chi_A(2n+1) \right)
= N - \sum_{n=0}^{2N-1} \chi_A(n)
\]
and so 
\[
\card(A \cap [0,2N-1]) = N
\]
which is condition~(i).
The second equation is equivalent to condition~(ii).  If this condition holds, then $Q(z)$ is a polynomial of degree at most $N-2.$  This completes the proof.  
\end{proof}

\bprob
Let $\ell \geq 3$  Does there exist a partition of the nonnegative integers into  pairwise disjoint sets $A_1, A_2,\ldots, A_{\ell}$ whose representation functions $r_{A_i,2}(n)$ for $i = 1, 2, \ldots, \ell$ eventually coincide?
\eprob

\section{Representation functions for sets of integers}

\subsection{Unique representation bases for the integers}
Sumsets of integers are very different from sumsets of nonnegative integers.  For example, the Erd\H os-Tur\' an conjecture asserts that the representation function of a basis of order 2 for the nonnegative integers must be unbounded.  In sharp contrast to this, there exist bases for the integers whose representation functions are bounded.  Indeed, we shall construct a basis $A$ of order 2 for \Z\ whose representation function is identically equal to 1.  Such sets are called \emph{unique representation bases}.

\bt[Nathanson~\cite{nath03a}]   \label{RS:theorem:URB}
Let $\varphi(x)$ be a function such that $\lim_{x\rightarrow\infty} \varphi(x) = \infty.$
There exists an additive basis $A$ for the group \Z\ of integers such that
\[
r_{A,2}(n) = 1 \qquad\mbox{for all $n \in \Z$,}
\]
and
\[
A(-x,x) \leq \varphi(x)
\]
for all sufficiently large $x$.
\et

\bpf
We shall construct an ascending sequence of finite sets
$A_1 \subseteq A_2 \subseteq A_3 \subseteq \cdots$
such that, for all $k\in \N$ and $n\in \Z$, 
\[
|A_k| = 2k \qquad\text{and}\qquad r_{A_k}(n) \leq 1 
\]
and
\[
r_{A_{2k}}(n) = 1 \qquad\mbox{if $|n| \leq k$.}
\]
It follows that the infinite set
\[
A = \bigcup_{k=1}^{\infty} A_k
\]
is a unique representation basis for the integers.

We construct the sets $A_k$ by induction.
Let $A_1 = \{0, 1\}$.  We assume that for some $k \geq 1$ we have
constructed sets
\[
A_1 \subseteq A_2 \subseteq \cdots \subseteq A_k
\]
such that $|A_k| = 2k$ and
\[
r_{A_k}(n) \leq 1 \qquad\mbox{for all $n \in \Z$.}
\]
We define the integer
\[
d_k = \max\{|a| : a\in A_k\}.
\]
Then
\[
A_k \subseteq [-d_k,d_k]
\]
and
\[
2A_k \subseteq [-2d_k, 2d_k].
\]
If both numbers $d_k$ and $-d_k$ belong to the set $A_k$,
then, since $0 \in A_1 \subseteq A_k$ and $d_k \geq 1,$ 
we would have the following two representations of 0 in the sumset $2A_k$:
\[
0 = 0 + 0 = (-d_k)+ d_k.
\]
This is impossible, since $r_{A_k}(0) \leq 1$, hence only one of the
two integers $d_k$ and $-d_k$ belongs to the set $A_k$.
It follows that if $d_k \not\in A_k,$ then
\[
\{2d_k,2d_k-1\} \cap 2A_k = \emptyset,
\]
and if $-d_k \not\in A_k,$ then
\[
\{-2d_k,-(2d_k-1)\} \cap 2A_k = \emptyset.
\]
Select an integer $b_k$ such that 
\[
b_k = \min\{|b| : b\not\in 2A_k\}.
\]
Then
\[
1 \leq b_k \leq 2d_k-1.
\]
To construct the set $A_{k+1}$, we choose an integer $c_k$
such that
\[
c_k \geq d_k.
\]
If $b_k \not\in 2A_k$, let
\[
A_{k+1} = A_k \cup \{b_k + 3c_k, -3c_k\}.
\]
We have
\[
b_k = (b_k+3c_k) + (-3c_k) \in 2A_{k+1}.
\]
If $b_k \in 2A_k$, then $-b_k \not\in 2A_k$ and we let
\[
A_{k+1} = A_k \cup \{-(b_k + 3c_k), 3c_k\}.
\]
Again we have
\[
-b_k = -(b_k+3c_k) + 3c_k \in 2A_{k+1}.
\]
Since
\[
d_k < 3c_k < b_k + 3c_k,
\]
it follows that $|A_{k+1}| = |A_k|+2 = 2(k+1)$.
Moreover,
\[
d_{k+1} = \max\{|a| : a\in A_{k+1}\} = b_k + 3c_k.
\]

For example, since $A_1 = \{0,1\}$ and $2A_1 = \{0,1,2\}$, it follows
that $d_1 = b_1 = 1.$  Then $b_1\in 2A_1$ but $-1= -b_1\notin 2A_1.$  Choose an integer  $c_1 \geq 1$ and let 
\[
A_2 = \{-(1+3c_1),0,1,3c_1\}.
\]
Then
\[
2A_2 = \{-(2+6c_1), -(1+3c_1), -3c_1, -1, 0, 1, 2, 3c_1, 1+3c_1, 6c_1  \}
\]
and  $d_2 = 1+3c_1$ and $b_2 = 2$.  Moreover, $r_{A_2}(n)=1$ if $|n| \leq 1.$

Assume that $b_k \not\in 2A_k$, hence
$A_{k+1} = A_k \cup \{b_k + 3c_k, -3c_k\}$.
(The argument in the case $b_k \in 2A_k$ and $-b_k \not\in 2A_k$
is similar.)
The sumset $2A_{k+1}$ is the union
of the following four sets:
\[
2A_{k+1} = 2A_k
\cup \left( A_k + b_k + 3c_k \right)
\cup \left( A_k - 3c_k \right)
\cup \{b_k, 2b_k + 6c_k,-6c_k\}.
\]
We shall show that these sets are pairwise disjoint.
If $u \in 2A_k$, then
\[
-2c_k \leq -2d_k \leq u \leq 2d_k \leq 2c_k.
\]
Let $a \in A_k$ and $v  = a + b_k + 3c_k \in A_k + b_k + 3c_k$.  The inequalities
\[
-c_k \leq -d_k \leq a \leq d_k \leq c_k
\]
and
\[
1 \leq b_k \leq 2d_k-1 \leq 2c_k-1
\]
imply that
\[
2c_k + 1 \leq v \leq 6c_k - 1 < 2b_k + 6c_k.
\]
Similarly, if $w = a - 3c_k \in A_k - 3c_k$, then
\[
-6c_k < -4c_k \leq w \leq -2c_k.
\]
These inequalities imply that the sets $2A_k$,
$A_k + b_k + 3c_k$, $A_k - 3c_k$, and $2\{b_k + 3c_k,-3c_k\}$
are pairwise disjoint, unless $c_k = d_k$ and
$-2d_k \in 2A_k \cap (A_k - 3d_k)$.
If $-2d_k \in 2A_k$, then $-d_k \in A_k$.
If $-2d_k \in A_k-3d_k$, then $d_k \in A_k.$
This is impossible, however, because
the set $A_k$ does not contain both integers $d_k$ and $-d_k$.

Since the sets $A_k + b_k + 3c_k$ and $A_k - 3c_k$ are
translations, it follows that
\[
r_{A_{k+1}}(n) \leq 1 \qquad\mbox{for all integers $n$}.
\]

Let $A = \bigcup_{k=1}^{\infty} A_k$.
For all $k \geq 1$ we have $2 = b_2 \leq b_3 \leq \cdots$
and $b_k < b_{k+2}$, hence $b_{2k} \geq k+1$.
Since $b_{2k}$ is the minimum of the absolute values of the integers
that do not belong to $2A_{2k}$, it follows that
\[
\{-k,-k+1,\ldots,-1,0,1,\ldots, k-1, k\} \subseteq 2A_{2k} \subseteq 2A
\]
for all $k \geq 1$, and so $A$ is an additive basis of order 2.   In particular, $r_{A_{2k}}(n) \geq1$ for all $n$ such that $|n| \leq k$.
If $r_{A,2}(n) \geq 2$ for some $n$,
then $r_{A_k,2}(n) \geq 2$ for some $k$, which is impossible.
Therefore, $A$ is a unique representation basis for the integers.

We observe that if $x \geq 1$ and $k$ is the unique integer such that 
$d_k \leq x < d_{k+1},$ then 
\bq
A(-x,x) 
& = & A_{k+1}(-x,x) \\
& = & \left\{
\begin{array}{ll}
2k & \mbox{for $d_k\leq x < 3c_k$,}  \\
2k+1 & \mbox{for $3c_k\leq x < b_k+3c_k = d_{k+1}$.}
\end{array}
\right.
\eq
In the construction of the set $A_{k+1}$, the only constraint 
on the choice of the number $c_k$ was that $c_k \geq d_k$.
Given a function $\varphi(x)$ such that $\lim_{x\rightarrow\infty} \varphi(x)=\infty,$
we shall use induction to construct a sequence of integers
$\{c_k\}_{k=1}^{\infty}$
such that $A(-x,x) \leq \varphi(x)$ for all $x \geq c_1.$
We begin by choosing a positive integer $c_1$ such that
\[
\varphi(x) \geq 4 \qquad \mbox{for $x \geq c_1$.}
\]
Then
\[
A(-x,x)  \leq 4 \leq \varphi(x) \qquad \mbox{for $c_1 \leq x \leq d_2$.}
\]
Let $k \geq 2$, and suppose we have selected an integer $c_{k-1} \geq d_{k-1}$ such that
\[
\varphi(x) \geq 2k \qquad \mbox{for $x \geq c_{k-1}$}
\]
and
\[
A(-x,x) \leq \varphi(x) \qquad \mbox{for $c_1 \leq x \leq d_k$.}
\]
There exists an integer $c_k \geq d_k$ such that
\[
\varphi(x) \geq 2k+2 \qquad \mbox{for $x \geq c_k$.}
\]
Then
\[
A(-x,x) = 2k \leq \varphi(x) \qquad \mbox{for $d_k \leq x < 3c_k$}
\]
and
\[
A(-x,x) \leq 2k+2 \leq \varphi(x) \qquad \mbox{for $3c_k \leq x \leq d_{k+1}$,}
\]
hence
\[
A(-x,x) \leq \varphi(x) \qquad \mbox{for $c_1 \leq x \leq d_{k+1}$.}
\]
It follows that
\[
A(-x,x) \leq \varphi(x) \qquad \mbox{for all $x \geq c_1$.}
\]
This completes the proof.
\epf

Theorem~\ref{RS:theorem:URB} constructs arbitrarily sparse unique representation bases.   If $A$ is a unique representation basis of order 2 with counting function $A(x),$ then $A(x) \ll x^{1/2}.$   We do not know how dense a unique representation basis can be.

\bprob
Let $\Theta$ be the set of all positive numbers $\theta$ such that there exists a unique representation basis $A$ with $A(x) \gg x^{\theta}.$  Compute $\sup \Theta.$
\eprob

There is work related to this problem by Chen~\cite{chen07} and Lee~\cite{lee07a}.
for all $x\in hA_1\setminus \{u_1\}.$

\subsection{Asymptotic bases for the integers}

Let $\mathcal{F}(\Z)$ denote the set of all functions from  \Z\ into $\N_0 \cup \{\infty\} .$   We shall consider the following two subsets of this function space:  The set of functions with only finitely many zeros, 
\[
\mathcal{F}_0(\Z)  = \{ f \in \mathcal{F}(\Z) : \card\left( f^{-1}(0) \right) < \infty \}
\]
and the set of functions that are nonzero for almost all integers $n$, 
\[
\mathcal{F}_{\infty}(\Z)  = \{ f \in \mathcal{F}(\Z) : d\left( f^{-1}(0) \right) =0 \}.
\]
For every positive integer $h,$ let $\mathcal{R}_h(\Z)$ denote the set of all representation functions of $h$-fold sumsets, that is, 
\[
\mathcal{R}_h(\Z) = \{ f\in \mathcal{F}(\Z) : f = r_{A,h} \text{ for some $A \subseteq \Z$} \}.
\]
For example, $\mathcal{R}_1(\Z) = \{ f: \Z \rightarrow \{0,1\} \}.$ 

Let $h \geq 2.$  If $A$ is a set of integers and $a \in A,$ then $r_{A,h}(ha) \geq 1$.   It follows that if $f \in \mathcal{F}(\Z)$ is a nonzero function such that $f(n) = 0$ for all $n \equiv 0 \pmod{h},$ then $f$ is not a representation function, and so $\mathcal{F}(\Z) \neq \mathcal{R}_h(\Z)$.

\bprob
Let $h \geq 2.$  Find necessary and sufficient conditions for a function $f \in \mathcal{F}(\Z)$ to be the representation function of an $h$-fold sumset.  
\eprob

This is called the \emph{inverse problem for representation functions in additive number theory}.

The set $A$ is an asymptotic basis of order $h$ for the integers if all but finitely many integers can be represented as the sum of $h$ not necessarily distinct elements of $A$.  Equivalently,  $A$ is an asymptotic basis of order $h$ for \Z\ if the representation function $r_{A,h}$ is an element of the function space $\mathcal{F}_0(\Z)$.  
We define
\[
\mathcal{R}_{h,0}(\Z) = \{ f\in \mathcal{F}_0(\Z) : f = r_{A,h} \text{ for some $A \subseteq \Z$} \}.
\]
Thus, $\mathcal{R}_{h,0}(\Z)$ is the set of representation functions of asymptotic bases of order $h$ for \Z.  
We shall prove the following important result:  For every integer $h \geq 2,$ 
\vspace{0.3cm}
\begin{center}
\framebox{ \Large{$\mathcal{R}_{h,0}(\Z) =  \mathcal{F}_0(\Z)$. } }
\end{center}
\vspace{0.3cm}
This means that \emph{every} function $f:\Z \rightarrow \N_0\cup \{\infty\}$ with only finitely many zeros is the representation function for some asymptotic basis of order $h$ for the integers.   

The proof will use Sidon sets.
A subset $A$ of an additive abelian semigroup $\X$ is called a \emph{Sidon set of order $h$}  if every element in the sumset $hA$ has a unique representation (up to permutations of the summands) as a sum of $h$ elements of $\X$.  Equivalently, $A$ is a Sidon set if $r_{A,h}(x) \leq 1$ for all $x \in \X.$  Sidon sets of order $h$ are also called \emph{$B_h$-sets}.  For example, every two-element set $\{a,b\}$ of integers (or two-element subset $\{a,b\}$ of any torsion-free abelian semigroup) is a Sidon set of order $h$ for all positive integers $h$, since the $h$-fold sumset
\[
h\{a,b\} = \{(h-i)a+ib : i=0,1,\ldots,h\} = \{ha+i(b-a):i=0,1,\ldots,h\}
\]
is simply an arithmetic progression of length $h+1$ and difference $b-a.$    Note that if the set $A$ is a Sidon set of order $h$, then $A$ is also a Sidon set of order $h'$ for all $h' = 1,2,\ldots, h-1.$

The set $A$ will be called a \emph{generalized Sidon set of order $h$} if, 
for all pairs of positive integers $r, r'$ with $r\leq h$ and $r' \leq h$, and for all sequences $a_1,\ldots,a_r$ and $a'_1,\ldots,a'_{r'}$ of elements of $A$, we have
\[
a_1 + \cdots  + a_r = a'_1 + \cdots + a'_{r'} 
\]
if and only if $r=r'$ and $a'_i = a_{\sigma(i)}$ for some permutation $\sigma$ of $\{1,\ldots,r\}$ and all $i = 1,\ldots, r$.

Note that if $A$ is a Sidon set (resp. generalized Sidon set) of order $h$, then $A$ is also a Sidon set (resp. generalized Sidon set) of order $h'$ for all positive integers $h' < h.$

\bl    \label{lemma:h-key}
Let $h \geq 2$ and let $c$ and $u$ be integers such that $c > 2h|u|.$  Then 
\[
D_{c,u} = \{-c,(h-1)c +u\}
\]
is a generalized Sidon set of order $h,$ and $u \in hD_{c,u}.$  Moreover, 
\[
\min\left\{|x-y|:x,y \in \bigcup_{r=1}^h rD_{c,u} \text{ and } x\neq y \right\} \geq c/2.
\]
\el

\bpf
We have
\[
u = (h-1)(-c)+((h-1)c+u) \in hD_{c,u}.
\]

To show that $D_{c,u}$ is a generalized Sidon set, let $i,j,i',j'$ be nonnegative integers such that
\[
1 \leq i+j \leq i'+j' \leq h.
\]
We define
\[
\Delta = \left[ i(-c)+j((h-1)c+u) \right] - \left[ i'(-c)+j' ((h-1)c+u) \right].
\]
If $\Delta = 0,$ then 
\[
(j'-j)hc = ((i'+j')-(i+j))c + (j-j')u.
\]
If $j'\neq j,$ then 
\begin{align*}
hc & \leq |(j'-j)hc| \\
& = | ((i'+j')-(i+j))|c + |j-j'||u| \\
& \leq (h-1)c+ h|u|  \\
& < \left(h - \frac{1}{2}\right) c
\end{align*}
which is absurd.  Therefore, $j=j'$ and so $i=i'$ and $D_{c,u}$ is a generalized Sidon set of order $h.$  

Suppose that $\Delta \neq 0.$  We must show that $|\Delta| > c/2.$  
If $j=j',$ then $i\neq i'$ and 
\[
|\Delta| = |i'-i|c \geq c.
\]
If $j \neq j',$ then 
\begin{align*}
|\Delta| & = \left| (j-j')hc  +  ((i'+j')-(i+j))|c + (j-j')|u\right|  \\
& \geq \left| j-j'\right|hc - \left|((i'+j')-(i+j) \right| c- \left|  (j-j')u\right|  \\
& \geq hc - (h-1)c - h |u| \\
& > \frac{c}{2}.
\end{align*}
This completes the proof.
\epf

\bt[Nathanson~\cite{nath04a,nath05a}]  \label{RS:theorem:FundRep}
Let $f:\Z \rightarrow \NO \cup \{ \infty\}$  be a function such that 
$\card\left( f^{-1}(0)\right) < \infty.$  For every $h \geq 2,$ there exists a set $A$ of integers such that $r_{A,h}(n) = f(n)$ for all $n \in \Z$.
\et

\bpf
We shall construct a sequence $\{A_k\}_{k=1}^{\infty}$ of finite sets  such that $A_k$ is a generalized Sidon set of order $h-1$ for all $k \geq 1,$ and $A = \cup_{k=1}^{\infty} A_k$ is an asymptotic basis of order $h$ for \Z\ whose representation function is equal to $f$.  

Let $U = \{u_k\}_{k=1}^{\infty}$ be a sequence of integers such that 
\[
\card\left(  \left\{ k \in \N: u_k = n \right\} \right) = f(n)
\]
for all integers $n$.  
It suffices to construct finite sets $A_k$ such that, for all integers $n,$ we have 
\beq  \label{RS:cond1}
r_{A_k,h}(n) \leq f(n)
\eeq
and
\beq  \label{RS:cond2}
r_{A_k,h}(n) \geq \card\left(  \left\{ i \in \{1,2,\ldots,k\} : u_i=n \right\} \right).
\eeq

Choose positive integers $d_1$ and $c_1$ such that
\[
 f^{-1}(0) \subseteq [-d_1, d_1]
\]
and 
\[
c_1 > 2h(  d_1 + |u_1| ).
\]
By Lemma~\ref{lemma:h-key}, the set 
\[
A_1 = D_{c_1, u_1} = \{-c_1, (h-1)c_1+u_1 \}
\]
is a generalized Sidon set of order $h$ and $u_1\in hA_1.$
We shall prove that $hA_1\cap f^{-1}(0) = \emptyset.$
If $x\in f^{-1}(0),$ then $|x| \leq d_1$ and so $|x-u_1| \leq d_1+|u_1|.$
Again by Lemma~\ref{lemma:h-key}, if $x\in hA_1\setminus \{u_1\},$ then
\[
|x-u_1| > \frac{c_1}{2} > h(d_1 + |u_1|) \geq 2(d_1 + |u_1|).
\]
It follows that $hA_1\cap f^{-1}(0) = \emptyset$, and so $r_{A_1,h}(n) \leq 1 \leq f(n)$ for all $n \in hA_1$ and $r_{A_1,h}(u_1) =1.$  
Thus, the set $A_1$ satisfies conditions~\eqref{RS:cond1} and~\eqref{RS:cond2}.

Let $k \geq 2,$ and assume that we have constructed a generalized Sidon set $A_{k-1}$ of order $h-1$ that satisfies conditions~\eqref{RS:cond1} and~\eqref{RS:cond2}.
Choose positive integers $d_k$ and $c_k$ such that 
\[
f^{-1}(0) \cup \bigcup_{r=1}^h rA_{k-1} \subseteq [-d_k,d_k]
\]
and
\[
c_k > 2h(2 d_k + |u_k| ).
\]
Let 
\[
A_k  = A_{k-1} \cup D_{c_k,u_k} = A_{k-1} \cup \{-c_k, (h-1)c_k + u_k\}.
\]
Then
\[
hA_k = hA_{k-1} \cup \bigcup_{r=1}^{h} \left(rD_{c_k,u_k} + (h-r)A_{k-1}\right).
\]
By Lemma~\ref{lemma:h-key}, the set $D_{c_k,u_k}$ is a generalized Sidon set of order $h$, and so every integer in the set $ \bigcup_{r=1}^{h} rD_{c_k,u_k}$ has exactly one representation as the sum of at most $h$ elements of $D_{c_k,u_k}$.  Also, the minimum distance between the elements of $ \bigcup_{r=1}^{h} rD_{c_k,u_k}$ is greater than $c_k/2.$  

Let $x,x' \in  \bigcup_{r=1}^{h} rD_{c_k,u_k}$ with $x\neq x'.$  
By Lemma~\ref{lemma:h-key}, there are unique positive integers $r,r'$ such that 
$x \in rD_{c_k,u_k}$ and $x' \in r' D_{c_k,u_k}$.  If $y \in (h-r)A_{k-1}$ and $y' \in (h-r')A_{k-1}$, then 
\[
|y-y'| \leq |y|+|y'| \leq 2d_k < \frac{c_k}{2} \leq |x'-x|
\]
and so $x+y \neq x'+y'$.  It follows that the sets $\{x\} + (h-r)A_{k-1}$ and $ \{x'\} + (h-r')A_{k-1}$ are pairwise disjoint.  Since $A_{k-1}$ is a generalized Sidon set of order $h-1$, it follows that every element of 
\[
\bigcup_{r=1}^{h} \left(rD_{c_k,u_k} + (h-r)A_{k-1}\right)
\]
has a unique representation as the sum of exactly $h$ elements of $A_k.$  

Recall that $u_k \in h D_{c_k,u_k}$ and 
$hA_{k-1} \cup f^{-1}(0) \subseteq [-d_k,d_k].$
If $w \in  hA_{k-1} \cup f^{-1}(0),$ then $|u_k-w| \leq d_k+|u_k|.$
If
\[
z \in \bigcup_{r=1}^{h} \left(rD_{c_k,u_k} + (h-r)A_{k-1}\right)
\]
then $z=x+y,$ where $x \in rD_{c_k,u_k}$ for some $r \in [1,h]$ and $y \in (h-r)A_{k-1}.$  If $z \neq u_k,$ then $x\neq u_k$.  It follows again from Lemma~\ref{lemma:h-key} that $|x-u_k| \geq c_k/2$ and 
\begin{align*}
|z-w| & = |x+y-w| = |x-u_k +u_k + y - w| \\
& \geq |x-u_k|-|u_k + y - w| \\
& \geq \frac{c_k}{2} - (2 d_k + |u_k| ) \\
& > (h-1)(2 d_k + |u_k| ) \\
& > 0.
\end{align*}
Therefore,
\[
hA_k \subseteq \Z\setminus f^{-1}(0)
\]
and 
\[
hA_{k-1} \cap \left( \bigcup_{r=1}^{h} \left(rD_{c_k,u_k} + (h-r)A_{k-1}\right)\right) 
= \emptyset \text{ or } \{u_k\}.
\]
It follows that
\[
r_{A_k,h}(n) = 
\begin{cases}
r_{A_{k-1},h}(n) & \text{if $n \in hA_{k-1} \setminus \{ u_k \}$} \\
r_{A_{k-1},h}(u_k) + 1 & \text{if $n = u_k $} \\
1 & \text{if $n \in hA_k \setminus hA_{k-1}$}
\end{cases}
\]
and so the set $A_k$ satisfies conditions~\eqref{RS:cond1} and~\eqref{RS:cond2}.

A similar argument shows that $A_k$ is a generalized Sidon set of order $h-1.$  Let
\begin{align*}
Z & = \bigcup_{h'=1}^{h-1} h'A_k 
= \bigcup_{h'=1}^{h-1} \left( \bigcup_{\substack{r,s=0\\r+s=h'}}^{h'} \left(rD_{c_k,u_k} + sA_{k-1}\right)\right) \\
& = \bigcup_{\substack{r,s=0\\1\leq r+s\leq h-1}}^{h-1} \left(rD_{c_k,u_k} + sA_{k-1}\right)
\end{align*}
Suppose that
\[
z = x+y = x'+y' \in Z
\]
where $x\in rD_{c_k,u_k}, y \in sA_{k-1}, x' \in r'D_{c_k,u_k}, y' \in s'A_{k-1}$ for nonnegative integers $r,s,r',s'$ such that $1 \leq r+s \leq r'+s' \leq h-1.$  If $x\neq x',$ then 
\[
|x-x'| \geq \frac{c_k}{2} > 2d_k \geq |y'-y|
\]
and so $x-x'\neq y'-y,$ which is absurd.  Therefore, x=x' and $y=y'.$  Since $D_{c_k,u_k}$ is a generalized Sidon set of order $h$ and $A_{k-1}$ is a generalized Sidon set of order $h-1,$ it follows that $x$ and $y$ have \emph{unique} representations as sums of at most $h-1$ elements of $D_{c_k,u_k}$ and $A_{k-1}$, respectively, and so $z$ has a unique representation as the sum of at most $h-1$ elements of $A_k.$  This completes the proof.
\epf

By Theorem~\ref{RS:theorem:FundRep}, for every function $f \in \mathcal{F}_0(\Z),$ there exist infinitely many asymptotic bases  $A$ of order $h$ such that $r_{A,h} = f,$
and such bases can be constructed that are arbitrarily sparse.  An open problem is to determine how dense such a set can be.  
Nathanson and Cilleruelo~\cite{nath07f,nath07h} proved that for every $f \in  \mathcal{F}_0$ and every $\varepsilon > 0,$ there is a set $A$ of integers with $r_{A,h} = f$ and 
\[
A(-x,x) \gg x^{\sqrt{2}-1-\varepsilon}
\]
for all $x \geq 1.$
The construction uses dense Sidon sets.

\bprob
Let $\alpha_2^{\ast}$ be the supremum of the set of all positive real numbers $\alpha$ such that, for every $f \in  \mathcal{F}_0,$ there is a set $A$ of integers with $r_{A,h} = f$ and $A(-x,x) \gg x^{\alpha}$ for all $x \geq 1.$  Determine  $\alpha_2^{\ast}$.
\eprob

\bprob
Let $h \geq 3.$  Does there exist a positive real number $\alpha_h$ such that, for every $f \in  \mathcal{F}_0,$ there is a set $A$ of integers with $r_{A,h} = f$ and $A(-x,x) \gg x^{\alpha_h}$ for all $x \geq 1.$  How large can $\alpha_h$ be?  
\eprob

We can extend the inverse problem for representation functions to functions $f:\Z \rightarrow \N_0 \cup \{\infty\}$  that have infinitely many zeros.  In the case $h=2,$ if $f^{-1}(0)$ is a set of integers of density 0, then there we can construct a set $A$ with  $f = r_{A,2}.$  The problem is open for higher orders $h$.

\bprob
Let $f:\Z \rightarrow \N_0 \cup \{\infty\}$ be a function such that 
$d\left( f^{-1}(0) \right) = 0.$  Let $h \geq 3.$  Does there exist a set $A$ of integers such that $r_{A,h}(n) = f(n)$ for all integers $n$?  
\eprob

We can extend this problem to functions whose zero sets have small positive density.

\bprob
Let $h \geq 2.$  Does there exist $\delta = \delta(h) > 0$ such that if $f:\Z \rightarrow \N_0 \cup \{\infty\}$  a function with 
$d_U\left( f^{-1}(0) \right) < \delta,$ then there exists a set $A$ of integers such that $r_{A,h}(n) = f(n)$ for all integers $n$?  
\eprob

\section{Representation functions for abelian semigroups}
The significant difference between inverse problems for \NO\ and \Z\ derives in part from the fact that \Z\ is a group but \NO\ is not.  Nathanson~\cite{nath04e} obtained some general inverse theorems for representation functions of ``semigroups with a group component.''

Let $B$ be a subset of an abelian semigroup \X\ and let $x\in \X.$  We define the representation functions 
\[
r_{B,2}(x) = \card\left( \left\{ \{b,b'\} \subseteq B : b+b'=x  \right\} \right)
\]
and
\[
\hat{r}_{B,2}(x) = \card\left( \left\{ \{b,b'\} \subseteq B : b+b'=x \text{ and } b \neq b' \right\} \right).
\]

We consider semigroups $S$ with the property that $S+S=S.$  
Equivalently, for every $s \in S$ there exist $s', s'' \in S$ such that
$s = s' + s''.$  Every semigroup with identity has this property,
since $s = s + 0.$  There are also semigroups without identity that have this property.  For example, if $S$ is any totally ordered set without a smallest element,
and if we define $s_1 + s_2 = \max(s_1,s_2),$ then $S$ is an abelian semigroup
such that $s = s+s$ for all $s \in S,$ but $S$ does not have an identity element.

Let $S$ be an abelian semigroup and let $B \subseteq S.$    For every positive integer $h$, we define the \emph{dilation}
\[
h \ast B = \{ hb : b\in B\}  = \{ \underbrace{b+ \cdots+b}_{\text{$h$ summands}} : b\in B\}.
\]
Note that if $G$ is an abelian group such that every element of $G$ has order dividing $h$, then $h\ast G = \{ 0\}.$

\bt  \label{rf:theorem:1}
Let $S$ be a countable abelian semigroup such that 
for every $s \in S$ there exist $s', s'' \in S$ with $s = s' + s''.$  
Let $G$ be a countably infinite abelian group such that 
the dilation $2\ast G$ is infinite.
Consider the abelian semigroup $\X = S \oplus G$  
with projection map $\pi:\X\rightarrow G$.
Let 
\[
f:\X \rightarrow \N_0 \cup \{\infty\}
\]
be any map such that the set 
$\pi\left(f^{-1}(0)\right)$
is a finite subset of $G$.  
Then there exists a set $B \subseteq \X$ such that 
\[
\hat{r}_{B,2}(x) = f(x)
\]
for all $x \in \X.$ 
\et

Note that Theorem~\ref{rf:theorem:1} is not true 
for all abelian semigroups.  For example, 
let \N\ be the additive semigroup of positive integers under addition,
and $\X = \N \oplus \Z.$  Since the equation $s' + s'' = 1$ has no solution in positive integers, it follows that, for every set $B \subseteq \X$, we have
$r_B(1,n) = \hat{r}_B(1,n) =  0$ for every $n \in \Z$.
Thus, if $f:\X \rightarrow \N_0 \cup \{\infty\}$ is any function with $f(1,n)\neq 0$ for some integer $n$, then there does not exist a set $B \subseteq \X$ with $\hat{r}_{B,2} = f.$

\bt  \label{rf:theorem:2}
Let $G$ be a countably infinite abelian group such that 
the dilation $2\ast G$ is infinite.
Let 
\[
f:G \rightarrow \N_0 \cup \{\infty\}
\]
be any map such that $ f^{-1}(0)$
is a finite subset of $G$.  
Then there exists a set $B$ of order 2 for $G$ such that 
\[
\hat{r}_{B,2}(x) = f(x)
\]
for all $x \in \X.$ 
\et

\bt\label{rf:theorem:3}
Let $S$ be a countable abelian semigroup such that 
for every $s \in S$ there exist $s', s'' \in S$ with $s = s' + s''.$  
Let $G$ be a countably infinite abelian group such that 
the dilation $12\ast G$ is infinite.
Consider the abelian semigroup $\X = S \oplus G$  
with projection map $\pi:\X\rightarrow G$.
Let 
\[
f:\X \rightarrow \N_0 \cup \{\infty\}
\]
be any map such that the set 
$ \pi\left(f^{-1}(0)\right)$
is finite. 
Then there exists a set $B \subseteq \X$ such that 
\[
r_{B,2}(x) = f(x)
\]
for all $x \in \X.$ 
\et

\bt  \label{rf:theorem:4}
Let $G$ be a countably infinite abelian group such that 
the dilation $12\ast G$ is infinite.
Let 
\[
f:G \rightarrow \N_0 \cup \{\infty\}
\]
be any map such that the set 
$f^{-1}(0)$
is finite.  
Then there exists an asymptotic basis $B$ of order 2 for $G$ such that 
\[
r_{B,2}(x) = f(x)
\]
for all $x \in \X.$ 
\et

The proofs of Theorems~\ref{rf:theorem:1}--\ref{rf:theorem:4}
can be found in~\cite{nath04e}.

\bprob
What countable abelian semigroups \X\ have the property that, for every function $f:\X \rightarrow \N_0 \cup \{\infty\}$ such that the set 
$f^{-1}(0)$ is finite, there exists an asymptotic basis $B$ of order 2 for $\X$ with $r_{B,2} = f$?
\eprob

\section{Bases associated to binary linear forms}

Let $\Phi(x_1,x_2) = u_1 x_1+u_2 x_2$ be a binary linear form with relatively prime integer coefficients $u_1$ and $u_2$.  Let $A_1$ and $A_2$ be sets of integers.  We define the set 
\[
\Phi(A_1,A_2) = \{\Phi(a_1,a_2) : a_1 \in A_1 \text{ and } a_2 \in A_2\}.
\]
The \emph{representation function} associated with the form ${\Phi}$ is 
\[
R_{A_1,A_2,{\Phi}}(n) = \card\left( \{ (a_1,a_2)\in A_1\times A_2 :  \Phi(a_1,a_2) = n \} \right).
\]
Then $R_{A_1,A_2,{\Phi}}$ is a function from \Z\ into $\N_0 \cup \{\infty\}.$
If $A_1=A_2=A,$  we write
\[
\Phi(A) = \Phi(A,A) = \{\Phi(a_1,a_2) : a_1,a_2 \in A\}
\]
and
\[
R_{A,{\Phi}}(n) = R_{A,A,{\Phi}}(n) 
= \card\left( \{ (a_1,a_2)\in A^2 :  \Phi(a_1,a_2) = n \} \right).
\]
The set $A$ will be a called a \emph{unique representation basis with respect to the form ${\Phi}$} if $R_{A,{\Phi}}(n) = 1$ for every integer $n$.

\bl  \label{group:lemma:fundamental}
Let $\Phi(x_1,x_2) = u_1 x_1+u_2 x_2$ be a binary linear form with relatively prime positive integer coefficients $u_1<  u_2.$  Let $A$ be a finite set of integers and let $b$ be an integer.  Then there exists a set $C$ with $A\subseteq C $ and $|C \setminus A| = 2$ such that
\beq  \label{group:rArC}
R_{C,{\Phi}}(n) = 
\begin{cases}
R_{A,{\Phi}}(b) + 1 & \text{if $n = b$} \\
R_{A,{\Phi}}(n)       & \text{if $n\in \Phi(A)\setminus \{b\}$} \\
1                   & \text{if $n\in \Phi(C)\setminus \left( \Phi(A) \cup \{b\}  \right)$}\\
0                   & \text{if $n\notin \Phi(C)$.}
\end{cases}
\eeq
\el

\bpf
Since $\gcd(u_1,u_2)=1,$ there exist integers $v_1$ and $v_2$ such that $\Phi(v_1,v_2) = u_1v_1+u_2v_2 = 1.$  Then 
\begin{align*}
\Phi(bv_1+u_2t,bv_2-u_1t) & = u_1(bv_1+u_2t)+u_2(bv_2-u_1t) \\
& = b(u_1v_1+u_2v_2) = b
\end{align*}
for all integers $t.$  Let $B =  \{ bv_1+u_2t,bv_2-u_1t \}.$
If $t\neq (b(v_2-v_1)/(u_1+u_2),$ then $bv_1+u_2t \neq bv_2-u_1t$ and $|B|=2.$  
We shall prove that there exist infinitely many integers $t$ such that $A \cap B = \emptyset$ and the set
$C = A \cup B$ 
satisfies conditions~\eqref{group:rArC}.

If $d = \max(\{ |a| : a\in A\}),$ then $|\Phi(a)| \leq (u_1+u_2)d$ for all $a\in A$.
The set $\Phi(C)$ is the union of the sets $\Phi(A), \Phi(A,B), \Phi(B,A),$ and $\Phi(B).$  

If $c \in \{ \Phi(a,bv_1+u_2t) : a\in A\},$ then there exists $a\in A$ such that
\[
c = u_1 a + u_2(bv_1+u_2t) = (u_1a + u_2v_1b) + u_2^2t 
\]
and so $c>d$ for all sufficiently large integers $t$. 

If $c' \in \{ \Phi(a, bv_2-u_1t ) : a \in A\},$ then there exists $a'\in A$ such that
\[
c' = u_1 a' + u_2(bv_2-u_1t ) = (u_1a' + u_2v_2b) - u_1 u_2 t 
\]
and so $c< -d$ for all sufficiently large integers $t$.  Therefore,
\[
\Phi(A) \cap \Phi(A,B) = \emptyset
\]
for all sufficiently large integers $t$.

For every integer $t$, the functions $\Phi(a,bv_1+u_2t)$ and $\Phi(a,bv_2-u_1t )$ are  strictly increasing functions of $a$.
Moreover, there exist $a,a' \in A$ such that 
$\Phi(a,bv_1+u_2t) = \Phi(a',bv_2-u_1t )$ if and only if
\[
(u_1a + u_2v_1b) + u_2^2t  = (u_1a' + u_2v_2b) - u_1 u_2 t 
\]
that is,  if and only if 
\[
(u_1+u_2)u_2t = u_1(a'-a) + u_2(v_2 - v_1)b
\]
and this identity holds only for finitely many $t$.
Thus, for all sufficiently large $t$ we have $R_{A,B,f}(n) \leq 1$ for all $n\in \Z.$

Similarly, if $c \in \{ \Phi(bv_1+u_2t,a) : a\in A\},$ then there exists $a\in A$ such that
\[
c = u_1 (bv_1+u_2t) + u_2 a = (u_1v_1 b + u_2 a) + u_1 u_2 t 
\]
and so $c>d$ for all sufficiently large integers $t$. 

If $c' \in \{ \Phi(bv_2-u_1t,a ) : a\in A\},$ then there exists $a'\in A$ such that
\[
c' = u_1 (bv_2-u_1t ) + u_2 a'= (u_1v_2b + u_2 a') - u_1^2 t 
\]
and so $c< -d$ for all sufficiently large integers $t$.  Therefore,
\[
\Phi(A) \cap \Phi(B,A) = \emptyset
\]
for all sufficiently large integers $t$.  
By the same method, we can prove that  for all sufficiently large $t$ 
we have $R_{B,A,f}(n) \leq 1$ for all $n\in \Z$ and 
\[
\Phi(A) \cap \Phi(B,A) = \Phi(A,B) \cap \Phi(B,A) = \emptyset
\]

Finally,  the set $\Phi(B) =  \{ \Phi(b',b'') : b',b'' \in B \}$ consists of the integers $b, (u_1v_2+u_2v_1)b+(u_2^2-u_1^2)t,(u_1+u_2)bv_1+u_2(u_1+u_2)t),$ and $(u_1+u_2)bv_2-u_1(u_1+u_2)t )$.  The coefficients of $t$ are the pairwise distinct integers $u_2^2-u_1^2, u_2(u_1+u_2),$ and $-u_1(u_1+u_2),$ and these are different from the numbers $-u_1u_2, -u_1^2, u_1u_2,$ and $u_2^2,$ which are the coefficients of $t$ in $\Phi(A,B)$ and $\Phi(B,A).$  It follows that $|\Phi(B)|=4$ and that the sets $\Phi(A), \Phi(B,A), \Phi(A,B),$ and $\Phi(B)\setminus \{ b\}$ are pairwise disjoint for all sufficiently large $t.$   This completes the proof.
\epf

\bt
Let $\Phi(x_1,x_2) = u_1 x_1+u_2 x_2$ be a binary linear form with relatively prime positive integer coefficients $u_1<  u_2.$   There exists a unique representation basis with respect to the form $\Phi,$ that is, a set $A$ of integers such that $R_{A,\Phi}(n)=1$ for all $n\in \Z.$
\et

\bpf
We shall construct an increasing sequence of finite sets $A_1\subseteq A_2 \subseteq \cdots$ such that $R_{A_k,f}(n) \leq 1$ for all $k\in \N$ and $n\in \Z,$ and $A = \bigcup_{k=1}^{\infty} A_k$ is a unique representation basis for $f.$  Let $A_1 = \{0,1\}.$  Then $\Phi(A_1) = \{0,u_1,u_2,u_1+u_2\}.$   Since $0 < u_1 < u_2 < u_1+u_2$, it follows that $|\Phi(A_1)| = 4$ and $R_{A_1,f}(n) \leq 1$ for all $n\in \Z.$ 

Let $A_k$ be a finite set of integers such that  $R_{A_k,\Phi}(n) \leq 1$ for all $n\in \Z.$   Let $b$ be an integer such that
\[
|b| = \min\left( \left\{ |n| : n \notin \Phi(A_k) \right\} \right).
\]
By Lemma~\ref{group:lemma:fundamental}, there is a set $A_{k+1}$ containing $A_k$ such that $b\in \Phi(A_{k+1})$ and $R_{A_{k+1},\Phi}(n) \leq 1$ for all $n\in \Z.$  This completes the proof.
\epf

More general results about representation functions of binary linear forms appear in Nathanson~\cite{nath07m}.

\bprob
Determine all $m$-ary linear forms $\Phi(x_1,\ldots,x_m) = u_1x_1+ \cdots + u_mx_m$ with nonzero, relatively prime integer coefficients such that there exists a unique representation basis with respect to $\Phi$.
\eprob

\bprob
Let $m \geq 2$ and let $\Phi(x_1,\ldots,x_m)$ be an  $m$-ary linear form with nonzero, relatively prime integer coefficients.  Let $f:\Z\rightarrow \N_0\cup \{\infty\}$ be a function such that $f^{-1}(0)$ is finite.  Does there exist a set $A$ such that $R_{A,{\Phi}} = f$?
\eprob

\bprob
Determine all $m$-ary linear forms $\Phi$ such that if $A$ and $B$ are sets of integers with $R_{A,{\Phi}} = R_{B,{\Phi}},$ then $A=B .$
\eprob

\bprob
Determine all $m$-ary linear forms $\Phi$ such that if $A$ and $B$ are finite sets of integers with $R_{A,{\Phi}} = R_{B,{\Phi}},$ then $A=B .$
\eprob

The last problem is related to work of Ewell, Fraenkel, Gordon,  Selfridge, and Straus~\cite{ewel68,gord-frae-stra62,self-stra58}.

\def\cprime{$'$} \def\cprime{$'$} \def\cprime{$'$}
\providecommand{\bysame}{\leavevmode\hbox to3em{\hrulefill}\thinspace}
\providecommand{\MR}{\relax\ifhmode\unskip\space\fi MR }
\providecommand{\MRhref}[2]{%
  \href{http://www.ams.org/mathscinet-getitem?mr=#1}{#2}
}
\providecommand{\href}[2]{#2}

\end{document}